\input amstex
\documentstyle{amsppt}
\magnification 1200
\vcorrection{-9mm}

\topmatter
\title     On alternating quasipositive links
\endtitle

\author    S.~Yu.~Orevkov
\endauthor

\address
Steklov Mathematical Institute, Gubkina 8, Moscow, Russia
\endaddress

\address
IMT, l'universit\'e Paul Sabatier, 118 route de Narbonne, Toulouse, France
\endaddress

\email
orevkov\@math.ups-tlse.fr
\endemail

\abstract
We prove that if a quasipositive link can be represented by
an alternating diagram satisfying the condition that no pair
of Seifert circles is connected by a single crossing, then
the diagram is positive and the link is strongly quasipositive.
\endabstract
\endtopmatter

\def\refBa     {1}
\def\refBBG    {2}
\def\refDEHL   {3}
\def\refDHL    {4}
\def\refDP     {5}
\def\refHay    {6}
\def\refMu     {7}
\def\refNa     {8}
\def\refOrGAFA {9}
\def\refOrQPCS {10}
\def\refRuTop  {11}
\def\refRuIII  {12}
\def\refRuPos  {13}
\def\refTr     {14}
\def\refYa     {15}

\def\thDHL   {1}
\def\thMain  {2}
\def\corMain {3}
\def\thTr    {4}
\def\remDEHL {5}
\def\thDEHL  {6}
\def\remFinal{7}
\def\quest   {8}

\def\eqYa   {1}
\def\eqMT   {2}
\def\eqDP   {3}
\def\eqDEHL {4}

\def\Null{\bold n}

\document

\head 1. Introduction \endhead

An $n$-braid is called {\it quasipositive} if it is a product of conjugates
of the standard generators $\sigma_1,\dots,\sigma_{n-1}$ of the
braid group $B_n$. A braid is called {\it strongly quasipositive} if it
is a product of braids of the form $\tau_{k,j}\sigma_j\tau_{k,j}^{-1}$
for $j\le k$ where $\tau_{k,j}=\sigma_{k}\sigma_{k-1}\dots\sigma_j$.
All links in this paper are assumed to be oriented links in the 3-sphere $S^3$.
A link is called {\it (strongly) quasipositive} if it is the braid closure of a
(strongly) quasipositive braid
(see [\refRuTop, \refRuIII]).

S.~Baader [\refBa, p.~268, Question (4)] asked:
{\sl Do quasipositive alternating links have positive diagrams?}
Note that positive diagrams represent strongly quasipositive links
(see [\refNa], [\refRuPos]) and alternating strongly quasipositive links
have positive alternating diagrams by [\refBBG, Cor.~7.3]. Notice also
that positive alternating diagrams are special (a diagram is called
{\it special} [\refMu] if its Seifert circles bound disjoint disks).

In this note we give an affirmative answer for a rather large class of alternating
links: those which have an alternating diagram whose number of Seifert circles
is equal to the braid index of the link.
We call such diagrams
{\it Diao--Hetyei--Liu\/} or {\it DHL diagrams\/} (and the corresponding links
{\it DHL links\/})
because these authors gave in [\refDHL] the following very nice and simple
characterization for them.

\proclaim{ Theorem \thDHL } {\rm([\refDHL, Thm.~1.1])}
An alternating diagram is DHL if and only if
there is no pair of Seifert circles connected by a single crossing.
\endproclaim

Our main result is the following.

\proclaim{ Theorem \thMain }  
Let $D$ be a DHL diagram of a quasipositive link. Then $D$ is positive.
\endproclaim

The proof is an easy combination of
results from [\refDP], [\refHay], [\refMu], [\refTr], and [\refYa] (see Section 2).
Theorems~\thDHL\ and \thMain\ allow to produce a lot of
examples of non-quasipositive links without any computations.

Since any positive diagram represents a strongly quasipositive link
(see [\refNa], [\refRuPos]), we obtain:

\proclaim{ Corollary \corMain } Let $L$ be a DHL link.
Then the following conditions are equivalent:
\roster
\item"(i)" $L$ is quasipositive;
\item"(ii)" $L$ is strongly quasipositive;
\item"(iii)" $L$ has a positive alternating diagram.
\endroster
\endproclaim

In Section 3 we generalize Theorem~\thMain\ to all alternating links whose
braid index is computed in [\refDEHL]; see Theorem~\thDEHL\ and Remark~\remFinal.

\smallskip\noindent{\bf Acknowledgement.}
I am grateful to Michel Boileau for useful discussions.

\head 2. Proof of the main theorem \endhead

Let $D$ be a connected link diagram. The {\it Seifert graph\/} of $D$ is the graph $G_D$
whose vertices correspond to Seifert circles and the edges 
correspond to the crossings. The sign of an edge is the sign of the corresponding
crossing.
A diagram $D$ is called {\it reduced} if $G_D$ does not have any edge
whose removal disconnects $G_D$.
Let $d(D)$ denote the sum of signs of all edges of a spanning tree of $G_D$,
and let $w(D)$ be the {\it writhe} of $D$, i.e., the sum of signs of all crossings.

For a link $L$, let $\sigma(L)$ and $\Null(L)$ be its signature and nullity
(the latter is the nullity of a symmetrized Seifert form on a {\sl connected}
Seifert surface).

\proclaim{ Theorem~\thTr } {\rm(Traczyk [\refTr])}
Let $D$ be a connected reduced alternating diagram of a link $L$. Then
$\sigma(L)=d(D)-w(D)$ and $\Null(L)=0$.
\endproclaim

\noindent
This formula for $\sigma(L)$ is given in [\refTr, Thm.~2(1)]
(the factor $1/2$ is erroneous there).
The fact that $\Null(L)=0$ (equivalently,
$\det(L)\ne0$) is proven [\refMu, Lem.~5.1] and in the appendix to [\refTr].
Otherwise it can be easily derived from [\refTr, Thm.~1].

\demo{ Proof of Theorem~\thMain }
Let $D$ be a DHL diagram of a quasipositive link $L$.
Then each connected component of $D$ is evidently a DHL diagram and
it represents a quasipositive link by [\refOrQPCS].
So, it is enough to consider the case when $D$ is connected.

Let $n$ be the braid index of $L$. By definition of DHL diagrams,
$D$ has $n$ Seifert circles. Hence, by [\refYa, Thm.~1] (see the discussion
of this theorem in the introduction to [\refYa]), $L$ can be represented by
an $n$-braid $\beta_1$ with
$$
    w(\beta_1)=w(D).      \eqno(\eqYa)
$$
By [\refHay, Thm.~1.2] $L$ can be represented
by a quasipositive $n$-braid $\beta_2$. Then Murasugi--Tristram inequality
[\refMu] for quasipositive braids can be reformulated as follows
(see [\refOrGAFA, Cor.~3.2])
$$
    1+\Null(L) \ge |\sigma(L)| + n - w(\beta_2).    \eqno(\eqMT)
$$
By Dynnikov--Prasolov Theorem [\refDP] (Generalized Jones Conjecture)
we have
$$
    w(\beta_1) = w(\beta_2)       \eqno(\eqDP)
$$
By combining (\eqYa)\,--\,(\eqDP) with Theorem~\thTr\ (note that any DHL diagram
is reduced), we obtain $|d(D)-w(D)|\le 1-n+w(D)$ whence $w(D)-d(D)\le 1-n+w(D)$,
i.e., $d(D)\ge n-1$. Recall that $d(D)$ is the sum of signs of all edges of
a spanning tree of $G_D$. Any spanning tree of $G_D$ has $n-1$ edges, hence
all its edges are positive. Since each edge of $G_D$ belongs to some spanning tree,
we conclude that all crossings of $D$ are positive.
Theorem~\thMain\ is proven.
\qed\enddemo

\head 3. A generalization of the main theorem \endhead

Let $D$ be an alternating diagram of a link $L$. Let $b=b(L)$
be the braid index of $L$ and $s=s(D)$ be the number of Seifert circles of $D$.
Define $d^{\pm}=d^{\pm}(D)$ as the number of edges of this sign
in a spanning tree of $G_D$, thus $d=d(D)=d^+-d^-$.

Let $\beta$ be a braid with $b$ strands realizing $L$.
Due to Dynnikov -- Prasolov Theorem [\refDP],
$w(\beta)$ does not depend on the choice of $\beta$, which allows us
to define the numbers $r^{\pm}=r^{\pm}(D)$ from the system of equations
$$
     r^+ + r^- = s - b,\qquad r^+ - r^- = w(D) - w(\beta).
$$
{\bf Remark \remDEHL.} The definition of the numbers
$r^{\pm}$ in [\refDEHL] is not quite clear but in all cases when they are
computed in [\refDEHL], they satisfy our definition; cf. [\refDEHL, Rem.~3.1--3.3].

\medskip
If $D$ is a DHL diagram, then $r^+ = r^- = 0$
(recall that in this case $w(D)=w(\beta)$ by [\refYa, Thm.~1]),
thus the following statement is a generalization of Theorem~\thMain.

\proclaim{ Theorem \thDEHL }
Let $D$ be a reduced alternating diagram of a quasipositive link $L$, and
$$
        2r^-(D) \le d^-(D).                     \eqno(\eqDEHL)
$$
Then $D$ is positive
{\rm(and hence $L$ is strongly quasipositive by [\refNa, \refRuPos])}.
\endproclaim

\demo{ Proof } Since the arguments are almost the same as for Theorem~\thMain,
we just write down the final computation. So, we have
$w(D)-d\le|\sigma|\le 1-b+w(\beta)$, hence
$$
  d+1 \ge w(D) - w(\beta) + b
  = (r^+ - r^-) + s - (r^+ + r^-) = s-2r^- \ge s - d^-
$$
whence $d^+\ge s-1$ and the result follows.
\qed\enddemo

\noindent
{\bf Remark \remFinal.} In all cases when
the braid index of a reduced alternating diagram is computed in [\refDEHL],
the inequality (\eqDEHL) holds,
in particular it holds for minimal diagrams of two-bridge links
and of alternating Montesinos links.

\medskip\noindent
{\bf Question \quest.} Does (\eqDEHL) hold for any reduced alternating diagram?


\Refs

\ref\no\refBa\by S.~Baader
\paper Slice and Gordian numbers of track knots
\jour Osaka J. Math. \vol 42 \yr 2005 \pages 257--271 \endref

\ref\no\refBBG\by M.~Boileau, S.~Boyer, C.~M.~Gordon
\paper Branched covers of quasi-positive links and L-spaces
\jour J. of Topology \vol 12 \yr 2019 \pages 536--576 \endref

\ref\no\refDEHL\by Y.~Diao, G.~Ernst, G.~Hetyei, P,~Liu
\paper A diagrammatic approach for determining the braid index of alternating links
\jour Arxiv:1901.09778 \endref

\ref\no\refDHL\by Y.~Diao, G.~Hetyei, P,~Liu
\paper The braid index of reduced alternating links
\jour Arxiv:1701.07366 \endref

\ref\no\refDP\by I.~A.~Dynnikov, M.~V.~Prasolov
\paper Bypasses for rectangular diagrams.
      A proof of the Jones conjecture and related questions
\jour Trudy mosk. mat. obshch. \vol 74 \yr 2013 \issue 1 \pages 115-173
\lang Russian \transl English transl.
\jour Trans. Moscow Math. Soc. \vol 74 \yr 2013 \pages 97--144 \endref

\ref\no\refHay\by K.~Hayden
\paper Minimal braid representatives of quasipositive links
\jour Pac. J. Math. \vol 295 \yr 2018 \pages 421--427 \endref

\ref\no\refMu\by K.~Murasugi
\paper On certain numerical invariant of link types
\jour Trans. Amer. Math. Soc. \vol 117
\yr 1965 \pages 387--422\endref

\ref\no\refNa\by T.~Nakamura
\paper Four-genus and unknotting number of positive knots and links
\jour Osaka J. Math. \vol 37 \yr 2000 \pages 441--451 \endref

\ref\no\refOrGAFA\by S.~Yu.~Orevkov
\paper Classification of flexible $M$-curves of degree $8$ up to isotopy
\jour GAFA -- Geom. Funct. Anal. \vol 12 \yr 2002 \pages 723--755 \endref

\ref\no\refOrQPCS\by S.~Yu.~Orevkov
\paper Quasipositive links and connected sums
\jour Funk. anal. i ego prilozh. \vol 54 \yr 2020 \issue 1 \pages 81--86
\lang Russian \transl English transl.
\jour Funct. Anal. Appl.  \toappear \endref

\ref\no\refRuTop\by L.~Rudolph
\paper Algebraic functions and closed braids
\jour Topology \vol 22 \yr 1983 \pages 191--201 \endref

\ref\no\refRuIII\by L.~Rudolph
\paper A characterization of quasipositive Seifert surfaces
(constructions of quasipositive knots and links, III)
\jour Topology \vol 31 \yr 1992 \pages 231--237 \endref

\ref\no\refRuPos\by L.~Rudolph
\paper Positive links are strongly quasipositive
\inbook in: Proceedings of the Kirbyfest, Berkeley, CA, USA, June 22-26, 1998
\publ University of Warwick \publaddr Warwick, UK
\bookinfo Geom. Topol. Monogr. 2 \pages 555--562
 \yr 1999 \endref

\ref\no\refTr\by P.~Traczyk
\paper A combinatorial formula for the signature of alternating diagrams
\jour Fundamenta Math. \vol 184 \yr 2004 \pages 311--316 \endref

\ref\no\refYa\by S.~Yamada
\paper The minimal number of Seifert circles equals the braid index of a link
\jour Invent. Math. \vol 89 \yr 1987 \pages 347--356 \endref

\endRefs
\enddocument